\documentclass[10pt]{article}
\usepackage{amsfonts,amssymb, graphicx,a4,bm,bbm,amsmath,bbm}

\newcommand{\zeq}{\setcounter{equation}{0}}
\usepackage{color}
\newtheorem{theo}{Theorem}[section]

\definecolor{GreenYellow}{cmyk}{0.15,0,0.69,0}
\definecolor{Yellow}{cmyk}{0,0,1,0}
\definecolor{Goldenrod}{cmyk}{0,0.10,0.84,0}
\definecolor{Dandelion}{cmyk}{0,0.29,0.84,0}
\definecolor{Apricot}{cmyk}{0,0.32,0.52,0}
\definecolor{Peach}{cmyk}{0,0.50,0.70,0}
\definecolor{Melon}{cmyk}{0,0.46,0.50,0}
\definecolor{YellowOrange}{cmyk}{0,0.42,1,0}
\definecolor{Orange}{cmyk}{0,0.61,0.87,0}
\definecolor{BurntOrange}{cmyk}{0,0.51,1,0}
\definecolor{Bittersweet}{cmyk}{0,0.75,1,0.24}
\definecolor{RedOrange}{cmyk}{0,0.77,0.87,0}
\definecolor{Mahogany}{cmyk}{0,0.85,0.87,0.35}
\definecolor{Maroon}{cmyk}{0,0.87,0.68,0.32}
\definecolor{BrickRed}{cmyk}{0,0.89,0.94,0.28}
\definecolor{Red}{cmyk}{0,1,1,0}
\definecolor{OrangeRed}{cmyk}{0,1,0.50,0}
\definecolor{RubineRed}{cmyk}{0,1,0.13,0}
\definecolor{WildStrawberry}{cmyk}{0,0.96,0.39,0}
\definecolor{Salmon}{cmyk}{0,0.53,0.38,0}
\definecolor{CarnationPink}{cmyk}{0,0.63,0,0}
\definecolor{Magenta}{cmyk}{0,1,0,0}
\definecolor{VioletRed}{cmyk}{0,0.81,0,0}
\definecolor{Rhodamine}{cmyk}{0,0.82,0,0}
\definecolor{Mulberry}{cmyk}{0.34,0.90,0,0.02}
\definecolor{RedViolet}{cmyk}{0.07,0.90,0,0.34}
\definecolor{Fuchsia}{cmyk}{0.47,0.91,0,0.08}
\definecolor{Lavender}{cmyk}{0,0.48,0,0}
\definecolor{Thistle}{cmyk}{0.12,0.59,0,0}
\definecolor{Orchid}{cmyk}{0.32,0.64,0,0}
\definecolor{DarkOrchid}{cmyk}{0.40,0.80,0.20,0}
\definecolor{Purple}{cmyk}{0.45,0.86,0,0}
\definecolor{Plum}{cmyk}{0.50,1,0,0}
\definecolor{Violet}{cmyk}{0.79,0.88,0,0}
\definecolor{RoyalPurple}{cmyk}{0.75,0.90,0,0}
\definecolor{BlueViolet}{cmyk}{0.86,0.91,0,0.04}
\definecolor{Periwinkle}{cmyk}{0.57,0.55,0,0}
\definecolor{CadetBlue}{cmyk}{0.62,0.57,0.23,0}
\definecolor{CornflowerBlue}{cmyk}{0.65,0.13,0,0}
\definecolor{MidnightBlue}{cmyk}{0.98,0.13,0,0.43}
\definecolor{NavyBlue}{cmyk}{0.94,0.54,0,0}
\definecolor{RoyalBlue}{cmyk}{1,0.50,0,0}
\definecolor{Blue}{cmyk}{1,1,0,0}
\definecolor{Cerulean}{cmyk}{0.94,0.11,0,0}
\definecolor{Cyan}{cmyk}{1,0,0,0}
\definecolor{ProcessBlue}{cmyk}{0.96,0,0,0}
\definecolor{SkyBlue}{cmyk}{0.62,0,0.12,0}
\definecolor{Turquoise}{cmyk}{0.85,0,0.20,0}
\definecolor{TealBlue}{cmyk}{0.86,0,0.34,0.02}
\definecolor{Aquamarine}{cmyk}{0.82,0,0.30,0}
\definecolor{BlueGreen}{cmyk}{0.85,0,0.33,0}
\definecolor{Emerald}{cmyk}{1,0,0.50,0}
\definecolor{JungleGreen}{cmyk}{0.99,0,0.52,0}
\definecolor{SeaGreen}{cmyk}{0.69,0,0.50,0}
\definecolor{Green}{cmyk}{1,0,1,0}
\definecolor{ForestGreen}{cmyk}{0.91,0,0.88,0.12}
\definecolor{PineGreen}{cmyk}{0.92,0,0.59,0.25}
\definecolor{LimeGreen}{cmyk}{0.50,0,1,0}
\definecolor{YellowGreen}{cmyk}{0.44,0,0.74,0}
\definecolor{SpringGreen}{cmyk}{0.26,0,0.76,0}
\definecolor{OliveGreen}{cmyk}{0.64,0,0.95,0.40}
\definecolor{RawSienna}{cmyk}{0,0.72,1,0.45}
\definecolor{Sepia}{cmyk}{0,0.83,1,0.70}
\definecolor{Brown}{cmyk}{0,0.81,1,0.60}
\definecolor{Tan}{cmyk}{0.14,0.42,0.56,0}
\definecolor{Gray}{cmyk}{0,0,0,0.50}
\definecolor{Black}{cmyk}{0,0,0,1}
\definecolor{White}{cmyk}{0,0,0,0}

\newcommand{\V}{\mathcal{V}}

\let\a=\alpha    
       
\let\m=\mu \let\n=\nu \let\o=\omega    \let\p=\pi 
 \let\s=\sigma \let\t=\tau 
  
\let\D=\Delta  \let\G=\Gamma  
\let\O=\Omega    
\let\Y=\Upsilon
\let\\=\noindent
\def\ee{\end{equation}}
\def\be{\begin{equation}}

\title{Perfect and separating Hash families: new bounds via the algorithmic cluster expansion local lemma}
\author{
Aldo Procacci and
Remy Sanchis\\
\\
\small{ Departamento de Matem\'atica UFMG}
\small{ 30161-970 - Belo Horizonte - MG - Brazil}
\\
\footnotesize{e-mails: aldo@mat.ufmg.br, rsanchis@mat.ufmg.br}
}
\date{}

\begin{document}
\maketitle

\begin{abstract}

We present new lower bounds for the size of perfect and
separating hash families ensuring their existence. Such new bounds
are based on the algorithmic cluster expansion improved version of
the Lov\' asz Local Lemma, which also implies that the Moser-Tardos
algorithm finds such hash families in polynomial time.

\end{abstract}

{\footnotesize

}
\vskip.5cm

\section{Introduction and results}\label{sec1}
In this initial section we will review rapidly the state of the art of
the Lov\'asz Local Lemma, a powerful tool in the framework of the probabilistic method in combinatorics, focusing specifically
on the recent
cluster expansion improvement of the  Moser-Tardos algorithmic version of the Lemma.
We then will recall the main results
in the literature concerning Perfect Hash Families and Separating Hash Families. Finally we will present the  results of the paper.

\subsection{Lov\'asz Local Lemma: state of the art}
The Lov\'asz Local Lemma (LLL) was
originally formulated by Erd\"{o}s and Lov\'asz in \cite{EL} and since then it has turned out  to be one of the most powerful tools in the framework of
the probabilistic method in combinatorics
to prove the existence of combinatorial objects with certain desirable properties.
The philosophy of the  Lemma is basically to consider  a collection of ``bad" events in some suitably defined probability space
whose occurrence, even of just one of them, prevents the existence of a certain  ``good" event (i.e. the combinatorial object under analysis).
Then  the Lemma
provides a sufficient condition which, once satisfied,  guarantees that there is
a strictly positive probability  that none of the bad  events occurs (so that the good event exists). Such sufficient condition can be  inferred
from the so-called {\it dependency graph} of the collection of events.
We remind that a dependency graph  for a
collection of random events $\cal ~B$
is a (simple and undirected) graph $G$  with vertex set $\cal B$ such that each
event $B\in \cal B$ is independent from the
$\s$-algebra generated by the collection of events ${\cal B} \setminus \Gamma^*_G(B)$
where $\Gamma^*_G(B)=\Gamma_G(B)\cup\{B\}$, with  $\Gamma_G(B)$ denoting the neighborhood of $B$ in $G$, i.e. the
set of vertices of $G$ which are connected
to the vertex $B$ by an edge of $G$.

\noindent
The connection between the LLL and the cluster expansion of the abstract polymer gas, implicitly
implied by an old paper by Shearer \cite{Sh}, has been
sharply pointed out in  \cite{SS} by Scott and Sokal who  also showed that the LLL (with dependency graph $G$)
can be viewed as  a reformulation of  the Dobrushin criterion \cite{D}
for the convergence of the cluster expansion of the hard-core lattice gas (on the same graph $G$).

\\In a later paper \cite{FP} Fern\'andez and Procacci improved the Dobrushin criterion
and this  has then been used  straightforwardly by Bissacot et al.
in \cite{BFPS} to obtain a correspondent  improved cluster expansion
version of the LLL
 (shortly CLLL).
Such new version of the LLL has been already implemented to
get  new bounds  on  several graph coloring problems (see  \cite{NPS} and \cite{BKP}).


\\As the original Lov\'asz Local Lemma by Erd\"os-Lov\'asz, the  improved cluster expansion version by Bissacot et al.
 given in \cite{BFPS} is ``non-constructive", in the sense that
it claims the existence of a certain event without explicitly exhibiting it. Nevertheless, an algorithmic version of the CLLL,
based on a  breakthrough paper by Moser and Tardos \cite{MT},  has been recently provided in \cite{Pe} and  \cite{AP}.

\subsubsection{Moser Tardos setting (general case)}\label{mtset}
In
the Moser Tardos setting all events in the   collection  $\cal B$ depend on
a finite family   $\cal V$ of mutually independent random variable with
$\O$ being  the sample space determined by these variables so that a outcome $\o\in \O$ is just a random evaluation
of all variables of the family $\cal V$.
Each  event $B\in \cal B$  is supposed to depend  only on some subset of the  variables  $\cal V$,
denoted by $vbl(B)$.
Since variables in $\V$ are assumed to be mutually independent, any two events $B,B'\in \cal B$
such that  $vbl(B)\cap vlb(B')=\emptyset$ are necessarily independent.\def\E{{\cal E}}
Therefore the family  $\cal B$ has a natural dependency graph, i.e.  the graph $G$
with vertex-set $\cal B$ and edge-set   constituted by the pairs
$\{B,B'\}\subset \cal B$  such that $vbl(B) \cap vbl(B') \neq \emptyset$.

\\In this setting Moser and Tardos  define the following random algorithm,
whose output, when (and if) it stops,  is  an evaluation of the variables of the family $\cal V$  (i.e. an outcome $\o\in \O$)
which avoids all  the events in the collection $\cal B$.
\vskip.2cm

\noindent
{\bf MT-Algorithm (gereral case)}.
%
%
%

\noindent- {\bf Step 0}:  Sample all random variables in the family $\cal V$.

\\Let $\o_0\in \O$ be the output.
\vskip.1cm
\\~~\,For  $k\ge 1$

\\- {\bf Step $\bm k$}:

a) Take $\o_{k-1}\in \O$ and check all   bad events in the family $\cal B$.

b) {\it  i)} ~If  some bad  event occurs, choose one, say $B$, and  resample  its  variables

~~~~~~~~$vbl(B)$
leaving unchanged the remaining variables.

\\~~~~~~~~{\it ii)}\, If no bad event occurs, stop the algorithm.

\\Let $\o_k\in \O$ be the output.

%
%
%
%

\vskip.2cm

\\We are now in a position to state the algorithmic version of the CLLL, which will be the basic tool to get our
 results on perfect and separating hash families. We remind that an independent set in a graph $G$  is a set of vertices of $G$ no two of which are connected by an edge of $G$.

\begin{theo}\label{t1}{\bf [Algorithmic  CLLL]}\label{PAP} Given a finite set $\cal V$ of mutually independent random variables, let $\cal  B$
be a finite set of events determined by these variables  with natural dependency graph $G$.
Let  $\bm{\mu}=\{\mu_B\}_{B\in \cal B}$ be a sequence of real numbers
 in $[0,+\infty)$. If, for each $B\in \cal B$,
$$
Prob(B)\;\leq\;  \frac{\mu_B}{\displaystyle\sum_{\substack{Y\subseteq\Gamma^*_G(B)\\Y\,{\rm independent\, in}\,G}} \prod_{B'\in Y}\mu_{B'}}
$$
then the MT-algorithm  reaches an assignment of values of the variables $\cal V$ such that none  of the events in $\cal B$ occurs.
 Moreover the expected total number of resampling  steps made by the MT-algorithm to reach this assignment is at most
$\sum_{B\in \cal B}\m_B$.
\end{theo}
The proof of Theorem \ref{t1} can be found in \cite{Pe} and \cite{AP}.

\subsection{Perfect Hash Families and Separating Hash Families}\label{hashsec}

Given a finite set $U$ we denote by $|U|$ its cardinality. Given an integer $k$, we denote shortly $[k]=\{1,2,\dots, k\}$.
 A collection of sets
$\{W_1,\dots, W_k\}$  such that $W_i\cap W_j=\emptyset$ for all $\{i,j\}\subset [k]$
will be called hereafter  a ``disjoint family".

\vskip.1cm
\\Let $n, w$ be integers such that $2\le w\le n$. We  denote by $P_w([n])$
the set of all subsets of $[n]$ with cardinality $w$.
\\Given $s, w_1,w_2,\dots,w_s$ integers such that $\sum_{i=1}^s w_i=w$,
we denote by  $P^*_{w}([n])$   the set whose elements
are the disjoint families  $S=\{W_1,\dots, W_s\}$ such that $W_i\subset [n]$ and $|W_i|=w_i$ for $i=1,\dots ,s$.

\vskip.1cm
\\Let $A$ be a $N\times n$ matrix.
Given  $W\in P_w([n])$ we denote by  $A|_W$  the $N\times w$ matrix formed by the $w$ columns of the matrix $A$
with indices in $W$. Analogously,
given a disjoint family $S=\{W_1,\dots, W_s\}\in P^*_w([n])$,  we denote by $A|_S$
the $N\times w$ matrix formed by the $w$ columns of the matrix $A$
with indices in $W_1\cup\dots\cup W_s$.

\vskip.2cm

\\\underline{\it Perfect hash family}. Let $X$  and $Y$ be  finite sets with  cardinality $|X|=n$ and $|Y|=m$. Let $w\in \mathbb{N}$ such that $2\le w\le n$.
Then a perfect hash family of size $N$ is a sequence $f_1,\dots, f_N$ of functions from $X$ to $Y$ such that for any subset
$W\subset X$ with cardinality $|W|=w$ there exists $i\in \{1,\dots, N\}$ such that $f_i$ is injective when restricted to $W$.
Such perfect hash family will be denoted by  $ {\rm PHF}(N; n,m,w)$.

\\A perfect hash family $ {\rm PHF}(N; n,m,w)$ is usually  viewed as a matrix $A$ with $N$ rows and $n$ columns,
with entries in the set of integers $[m]\equiv\{1,2,\dots,m\}$ such that
for any set  $W\in P_w([n])$, the $N\times w$ matrix $A|_W$ formed by the $w$ columns of the matrix $A$
with indices in $W$ has at least one line  with
distinct entries.

\vskip.1cm

\\\underline{\it Separating  hash family}. Given   $X$  and $Y$   finite sets with  cardinality $|X|=n$ and $|Y|=m$  and the integers
$w_1,\dots,w_s$ such that  $w=w_1+\dots +w_s\le n$, a separating hash family of size $N$ is a sequence
$f_1,\dots, f_N$ of functions from $X$ to $Y$ such that for all disjoint families of subsets $\{W_1,\dots, W_s\}$
of $X$ such that $|W_j|=w_j$ ($j=1, \dots, s$), there exists $i\in \{1,\dots, N\}$ such that $\{f_i(W_1),\dots, f_i(W_s)\}$
is a disjoint family  of subsets of  $Y$.

\\A separating hash family $ {\rm SHF}(N; n,m,\{w_1,\dots, w_s\})$ can be viewed as a
matrix $A$ with $N$ rows and $n$ columns, with entries in the set of integers $[m]$  such that
for any disjoint family $S=\{W_1,\dots, W_s\}\in P^*_{w}([n])$,
the $N\times w$ matrix $A|_S$ formed
by the $w$ columns of the matrix $A$ with indices in $W_1\cup\dots\cup W_s$
has at least one line which  ``separate
$A|_{W_1},\dots, A|_{W_s}$", i.e.,
for any unordered  pair $\{r,r'\}\subset [s]$,
the entries of this line belonging to $A|_{W_r}$ are different from the entries
of the same line belonging to  $A|_{W_{r'}}$.
\vskip.2cm
\\The probabilistic method has been already used several times in
the past to face the problem of  the existence of perfect and separated hash families.
In particular lower  bounds for  $N$, ensuring the existence of a perfect hash family with fixed values $n,m, w$
have been first obtained
by Mehlhorn in \cite{Me} using standard techniques of the probabilistic method.
The  Local Lov\'asz Lemma  has been subsequently used by Blackburn \cite{Bl} to improve the Mehlhorn bound.
In the same year, another technique in the framework of the probabilistic method, the so-called expurgation method, has been used to get
alternative bounds for perfect hash families \cite{STW}.
Later the Lov\'asz Local Lemma has also been used
in \cite{DSW}  to get similar bounds also for separating hash families. In the same paper \cite{DSW} the authors also outlined a comparison
between  the LLL and the expurgation method for perfect hash families suggesting that the expurgation method yields better bounds
than the LLL.
In a related paper \cite{SZ} an alternative technique still based on the expurgation method has been used to obtain  new
lower  bounds for  $N$, for fixed values $n, m, \{w_1,w_2\}$, guaranteeing   the existence of separating hash families. We finally mention that
there have been also several results regarding  upper bounds for $N$ ensuring the non-existence of Separating and Hash families
(see, e.g., \cite{BT} and references therein)

\subsection{Results}
\\We conclude this introductory section by presenting our main
results which consist in new  bounds for perfect hash families
and separating hash families.

\\Our  first result concerns a lower bound for perfect hash families.

\begin{theo}\label{hash}
Let $N,n,n$ be integers and let $w$  be integer such that $2\le w\le n$. Then
there exists a perfect hash family ${\rm PHF}(N; n,m,w)$ as soon as
\be\label{henne}
N\ge ~\frac{\ln[\varphi_{w,n}'(\t)]+(w-1)\ln{(n-w)}-\ln{(w-1)!}}{\ln (m^{w})-\ln\left(m^{w}- w!{\binom{m}{w}}\right)}
\ee
where $\t$ is the first positive solution of the equation $\varphi_{w,n}(x)-x\varphi_{w,n}'(x)=0$ and
\be\label{fix}
 \varphi_{w,n}(x)= 1+\sum_{k=1}^{\lfloor n/w\rfloor \wedge w} {\binom{w}{k}} \tilde \G_k(w,n)x^{k}
\ee
 with
 $$
\tilde\G_k(w,n) =
\sum_{j=0}^{w-k}{\binom{w-k}{j}}
\prod_{\ell=1}^{k(w-1)-j-1}\left(1-{\frac{\ell}{n-w}}\right)\left[\frac{w}{n-w}\right]^{j}\times
$$
\be\label{fixg}
\underline{}~~~\times
\sum_{\substack{i_1+\dots+i_k=j\\ i_s\geq 0}} {\frac{j!}{i_1!\dots i_k!}}\prod_{\ell=1}^k\left[{\frac{1}{i_l+1}}\prod_{s=1}^{i_\ell}\Big(1-{\frac{s}{w}}\Big)\right]
\ee
Moreover  the MT-algorithm (described in Sec. \ref{halg} below)
finds such perfect hash family ${\rm PHF}(N; n,m,w)$  in an expected time which is polynomial in the input parameters $N$, $n$ and $m$ for
any fixed $w$.
\end{theo}

\\The second result concerns a similar lower bound for separating  hash families. To state this result
 we need to introduce the following definition. Given a multi-set
$w_1,\dots, w_s$ of integers such that $w_1+\dots+w_s=w$,   we denote by $m_p$ the multiplicity of the integer $p\in \{1,2,\dots,w\}$
in the multi-set $w_1,\dots, w_s$, i.e.
$m_p=\sum_{i=1}^s\mathbbm{1}_{\{w_i=p\}}$.

\begin{theo}\label{shash}
Let $N,n,n$ be integers and let $w$  be an integer such that $2\le w\le n$. Let $s\ge 2$ and let $\{w_1,\dots,w_s\}$ be a family of integers such that
$w_1+\dots+w_s=w$. Then
there exists a separating hash family $ {\rm SHF}(N; n,m,\{w_1,\dots, w_s\})$  as soon as

\be\label{SHF}
N\ge {\frac{\ln[\varphi_{w,n}'(\t)]+(w-1)\ln{(n-w)}-\ln{(w-1)!}+ \ln (m_{w})}{\ln\left({\frac{1}{q}}\right)}}
\ee
where $\varphi_{w,n}'(\t)$ is the same number introduced  in Theorem \ref{hash},
\be\label{mw}
m_{w} = \frac{1}{\prod_{p=1}^w m_p!}\frac{w!}{ w_1!\cdots w_k!}
\ee
and
$$
q=  1- \frac{\p_{\mathcal{G}_s}(m)}{ m^w}
$$
with $\p_{\mathcal{G}_s}(m)$ being
the chromatic polynomial of the
complete $s$-partite  graph $\mathcal{G}_s$  with  $w_1,\dots w_s$ vertices.
Moreover  the MT-algorithm (described in Sec. \ref{halg} below) finds such separating hash family ${\rm SHF}(N; n,m,\{w_1,\dots, w_s\})$
in an expected time which is polynomial in the input parameters $N$, $n$ and $m$ for
any fixed $w$.
\end{theo}

\\As claimed  in the abstract, we will use
the algorithmic version of the CLLL, i.e., Theorem \ref{t1}, to prove Theorems \ref{hash} and \ref{shash}. Let us thus conclude this section
by describing  how to adapt the Moser-Tardos setting and the MT-algorithm to the case of Perfect and separated hash families.

\subsubsection{Moser Tardos setting for (perfect [separating] hash families)}\label{halg}

\\In the present case of Hash families.
the finite family $\cal V$  of the Moser-Tardos setting is  constituted by a set
of $N n$ mutually independent random variable  taking values
in the set $[m]$ according to the uniform distribution and  representing the possible entries of a $N\times n$ matrix.
The sample space generated by the family $\cal V$ is thus $\O=[m]^{N\times n}$ and an outcome in $\O$ is
a $N\times n$ matrix $A$.

\\{\it The bad events}. For each   $W\in P_w([n])$  {[for each
$S=\{W_1,\dots, W_s\}\in P^*_{w}([n])$]},   let $E_W$  be the event such that in any line of
$A|_W$ at least two entries are equal [let $E_S$  be the event  such that
for any line of $A|_S=A|_{\cup_{r=1}^s W_r}$ there is a pair $\{r,r'\}\subset [s]$
such that two entries of this line,
one in $A|_{W_r}$ and the other in $A|_{W_{r'}}$, are equal]. We have thus a
family ${\cal W}\equiv \{E_W\}_{W\in P_w([n])}$ [a  family ${\cal S}\equiv \{E_S\}_{S\in P^*_{w}([n])}$] of bad events  containing
${n\choose w}$ members [containing ${n\choose w}m_{w}$ members, with $m_{w}$ defined in
(\ref{mw})]. If  $A$ is  a sampled  matrix  such that no bad event of the family
${\cal W}$ [${\cal S}$] occurs, then
 for every $W\subset P_w([n])$  [for every $S\in
P^*_{w}([n])$] at least one line of $A|_W$ [of $A|_S$] has distinct entries [separates $S=\{W_1,\dots, W_s\}$],
that is to
say $A$ is a ${\rm PHF}(N;n,m,w)$ [${\rm SHF}(N; n,m,\{w_1,\dots, w_s\})$].

\\We are now in the position to outline the {MT-algorithm for ${\rm PHF}$} [${\rm SHF}$]

\vskip.15cm

\noindent
{\bf MT-Algorithm (for ${\rm PHF}$ [for ${\rm SHF}$])}.
%
%
%

\noindent- {\bf Step 0}: Pick an evaluation  all $Nn$  variables of the family $\cal V$ (the matrix entries)

\\Let $A_0$ be the output matrix.
\vskip.1cm
\\~~\,For  $k\ge 1$

\\- {\bf Step $\bm k$}:

a) Take the matrix $A_{k-1}$ and check all   bad events of  the family $\cal W$

~~~~[of the family $\cal S$].

b) {\it  i)} \,If  some bad  event occurs, choose one, say $E_W$ [$E_S$], and take a new

~~~~~~~random evaluation of the entries of $A_{k-1}|_W$ [of the entries of $A_{k-1}|_S$]

~~~~~~~leaving unchanged the remaining entries of $A_{k-1}$.

\\~~~~~~~~{\it ii)}\,\,If no bad event occurs, stop the algorithm.

\\Let $A_k$ be the output matrix.

\vskip.15cm

%
\\Note that when the algorithm stops  the output  matrix is
a ${\rm PHF}(N;n,m,w)$ [the output  matrix is
a ${\rm SHF}(N;n,m,w)$]. We will see in the next section
this algorithm  stops after an expected number of steps equal to
${n\choose w}$.

\vskip.2cm

\\The rest of the
paper is organized as follows. In Section \ref{sechash} we give  the
proofs of Theorems \ref{hash} and \ref{shash}. Finally, in Section
\ref{compar} we discuss some comparisons with previous bounds given
in the literature.

\section{Proofs of Theorems \ref{hash} and \ref{shash}}\label{sechash}
\zeq

\subsection{Proof of Theorem \ref{hash}}

\\Let us  apply  Theorem \ref{PAP} for the family of events ${\cal W} =\{E_W\}_{W\in P_w([n])}$ introduced in the previous section (Sec. \ref{halg}).
Clearly two events
$E_W, E_{W'}\in \cal W$ are independent if $W\cap W'=\emptyset$. Therefore
the dependency graph $G$ of the family of events $\{E_W\}_{W\in P_w([n])}$
can be identified with the graph whose vertices are the elements $W$ of $P_w([n])$, i.e. the subsets
$W$ of $[n]$ with cardinality $w$, and two vertices $W$ and $W'$ of $G$
are connected by an edge of the dependency graph $G$  if and only if $W\cap W'\neq\emptyset$.
Thus the neighborhood of $W$ in $G$ is the set
$$
\G_G(W)=\{{W'}: ~W'\in  P_w([n])~{\rm and}~ W'\cap W\neq \emptyset\}
$$
The probability of an event $E_W$ is
$$
P(E_W)~=~\frac{[m^{w}- m(m-1)\cdots (m-w+1)]^N}{ m^{wN}} ~= ~\frac{\left[m^{w}- w!\binom{m}{ w}\right]^N}{ m^{wN}}
$$
and, according to Theorem \ref{PAP},
the MT-algorithm (as it was described in Section \ref{halg}) finds a   perfect hash family $ {\rm PHF}(N; n,m,w)$ if,  for
some $\m>0$
\be\label{FP0}
P(E_W)~=~\frac{\left[m^{w}- w!\binom{m}{w}\right]^N}{  m^{wN}}~\le~ \frac{\m_W}{ \displaystyle\sum_{\substack{Y\subseteq\Gamma^*_G(W)\\
Y~{\rm independent\,in\,}G}} \prod_{W'\in Y}\m_{W'}}
\ee
We set $\m_W= \m$ for all $W\in P_w([n])$ so that
$$
{ \displaystyle\sum_{\substack{Y\subseteq\Gamma^*_G(W)\\
Y~{\rm independent\,in\,}G}} \prod_{W'\in Y}\m_{W'}}~
= ~\displaystyle\sum_{\substack{Y\subseteq\Gamma^*_G(W)\\ Y~{\rm independent\,in\,}G}} \m^{|Y|}~=~ 1+ \sum_{k=1}^{\lfloor n/w\rfloor \wedge w} \G_k(w,n) \mu^k
$$
where
\be\label{gk}
\G_k(w,n)=\sum_{\substack{Y\subset\Gamma^*_G(W):\;|Y|=k,\\ Y{\rm
independent \,in\,}G}}1
\ee
Hence (\ref{FP0}) rewrites
\be\label{FP}
P(E_W)~=~\frac{\left[m^{w}- w!\binom{m}{ w}\right]^N}{  m^{wN}}~\le ~
\frac{\m}{\displaystyle1+ \sum_{k=1}^{\lfloor n/w\rfloor \wedge w} \G_k(w,n)~\mu^k}
\ee
Let us now calculate explicitly the number $\G_k(w,n)$
defined in (\ref{gk}). We have:

\begin{align*}
\G_k(w,n)&=
\frac{1}{ k!}\sum_{\substack{i_0+i_1+\dots+i_k=w\\ i_s\geq 1,~ s\geq 1}}
\frac{w!}{i_0!i_1!\cdots i_k!}\binom{n-w}{w-i_1} \binom{n-\!w-\!(w-i_1)}{w-i_2}~\cdots\\
&\;\;\;\;\;\;\;\;\;\;\;\;\;\;~~~~~~~~~~~~~~~~~~~~~~
\cdots\!  \binom{n-\!w-\!(w-i_1)-\!\dots\!-\!(w-i_{k-1})}{w-i_k}\\
&= \frac{1}{k!}\left(\frac{1}{w!}\right)^{\!\!k}
\sum_{i_0=0}^{w-k}\frac{w!}{i_0!} \sum_{\substack{i_1+\dots+i_k=w-i_0\\i_s\geq 1}}  \prod_{l=1}^{k} \binom{ w}{ i_l}
\frac{(n-w)!}{(n-kw-i_0)!}
\\
&= \binom{w}{k}\left[\frac{(n-w)^{w-1}}{w!}\right]^{k}
\sum_{j=0}^{w-k}\binom{w-k}{j}
\frac{\prod_{\ell=1}^{k(w-1)-j-1}\left(1-\frac{\ell}{ n-w}\right)}{ (n-w)^{j}}~~\times\\
&~~~~~~~~~~~~~~~~~~~~~~~~~~~~~~~~~~~~\times~~ j!\sum_{\substack{i_1+\dots+i_k=j\\ i_s\geq 0}}
\prod_{l=1}^{k}\binom{ w}{ i_l+1}\\
&
=  \binom{w}{k}\left[\frac{(n-w)^{w-1}}{w!}\right]^{k}
\sum_{j=0}^{w-k}\binom{w-k}{j}
\frac{\prod_{\ell=1}^{k(w-1)-j-1}\left(1-\frac{\ell}{n-w}\right)}{ (n-w)^{j}}~~\times\\
&\;\;\;\;\;\;\;\;\;\;\;\;\;\;\times
\sum_{\substack{i_1+\dots+i_k=j\\i_s\geq 0}}\frac{j!}{ i_1!\cdots i_k!}
\prod_{l=1}^{k}\left[\frac{ w!}{  (i_l+1)(w-i_l-1)!}\right]\\
&= \binom{w}{k}\!\!\!\left[\frac{(n-w)^{w-1}}{(w-1)!}\right]^{k}
\sum_{j=0}^{w-k}\binom{w-k}{ j}\left[\frac{w}{ n-w}\right]^{j}~
{\prod_{\ell=1}^{k\!(w-1)\!-\!j\!-\!1}\!\!\!\!\left(1-\frac{\ell}{ n-w}\right)}
\times\\
&\;\;\;\;\;\;\;\;\;\;\;\;\;\;\times
\sum_{\substack{i_1+\dots+i_k=j\\ i_s\geq 0}}\frac{j!}{ i_1!\cdots i_k!}
\prod_{l=1}^{k}\left[\frac{\prod_{s=1}^{i_l}(1-\frac{s}{ w})}{  (i_l+1)}\right].
\end{align*}
I.e.
we get
\be\label{GGT}
\G_k(w,n)= \binom{w}{ k} \left[\frac{(n-w)^{w-1}}{(w-1)!}\right]^{k}  \tilde \G_k(w,n),
\ee
with $\tilde \G_k(w,n)$ given by (\ref{fixg}).
Therefore, setting
$
\a= \frac{(n-w)^{w-1}}{(w-1)!}\,\mu
$
the condition (\ref{FP}) becomes

\be\label{FP2}
\frac{(n-w)^{w-1}}{(w-1)!}\left[\frac{m^{w}- w!\binom{m}{w}}{  m^{w}}\right]^N\le
\max_{\a>0}\frac{\a}{ \varphi_{w,n}(\a)}= \frac{1}{ \varphi_{w,n}'(\t)}
\ee
where
$$
 \varphi_{w,n}(\a)= 1+\sum_{k=1}^{\lfloor n/w\rfloor \wedge w} \binom{w}{k} \tilde \G_k(w,n)\a^{k}
 $$
and $\t$ is the first positive solution of the equation $\varphi(x)-x\varphi'(x)=0$.
Taking the logarithm on both sides of  (\ref{FP2}), condition  (\ref{FP}) is   thus implied by the following inequality
\be\label{FP32}
\frac{}{}N \ge\frac{A_n(w)}{ D_m(w)}
\ee
where
\be\label{anw}
A_n(w)~=~\ln[\varphi_{w,n}'(\t)]+(w-1)\ln{(n-w)}-\ln{(w-1)!}
\ee
and
\be
D_m(w)~=~ \ln (m^{w})-\ln\left(m^{w}- w!\binom{m}{ w}\right)~~~~~~~~~~~~~~~~~~
\ee
In conclusion once  (\ref{FP32}) is satisfied then also (\ref{FP0}) is satisfied and therefore,
according to Theorem \ref{t1},    a ${\rm PHF}(N;n,m,w)$ exists and  the MT-algorithm
(as described in  Section \ref{halg}) finds it
in an expected number of steps
$$\sum_{W\in P_w([n])} \m= |P_w([n])|\mu = {n\choose w}\mu\le {n\choose w}$$
where the last inequality follows from the fact that the optimum $\mu$ which maximize
the r.h.s. of (\ref{FP}) is surely less than one.

\\Now,  at  each step $k\ge 1$ of the MT-algorithm described in section \ref{halg}, in order
to check in item a) whether or not a bad event of the family $\{E_W\}_{P_w([n]}$ occurs,
we need to consider   all the $N$  lines of (at worst)  all the ${n\choose w}$  matrices  $A|_W$ with $W\in P_w([n]$,  and for
each line of a given matrix $A|_W$ we need to compare all pair of entries of the line to check
whether  they  are equal or not. This is done  in (at most) $N{w\choose 2}{n\choose w}$ operations.

\\This concludes the proof of Theorem \ref{hash}.

\subsection{Proof of Theorem \ref{shash}}

We first recall that, for a fixed sequence of integers $w_1,w_2,\dots,w_s$ such that $w= w_1+\dots w_s\le n$,
$P^*_{w}([n])$ is the set whose elements are
the disjoint families   $S=\{W_1,\dots, W_s\}$ of subsets of $[n]$ with cardinality $w_1,\dots, w_s$ resp. and $A|_S$
is the $N\times w$ matrix formed
by the $w$ columns of the matrix $A$ with indices in $\bigcup_{i=1}^s W_i$.  Given two disjoint families
$S=\{W_1,\dots, W_k\}$ and $S'=\{W'_1,\dots, W'_k\}$, we also denote shortly $S\cap S'\doteq
(\bigcup_{i=1}^s W_i)\cap  (\bigcup_{i=1}^s W'_i)$.

\\Let us  apply  Theorem \ref{PAP} for the family of events ${\cal S} =\{E_S\}_{S\in P^*_w([n])}$
introduced in Section \ref{halg}.

\\For  $S=\{W_1,\dots, W_s\}\in P^*_{w}([n]) $,  the probability of the event $E_S$ is given by
$$
P(E_S) =\prod_{i=1}^N P_i(S)
$$
where $P_i(S)$ is the probability that the line $i$ of the matrix $A|_S$ do not separate $S$. To calculate  $P_i(S)$ just observe that

$$
P_i(S)=\frac{\#\; \mbox{favorable cases}}{ \# \;\mbox{number of all cases}} =  1 -\frac{\#\; \mbox{unfavorable cases}}{ \# \;\mbox{number of all cases}} =
$$
Setting $w=w_1+w_2+\dots +w_s$ we have that
$$
\# \;\mbox{number of all cases}= m^w
$$
To count the number of unfavorable cases, consider the complete $s$-partite  graph
$\mathcal{G}_s$  with vertex set $W$ and independent sets of  vertices set $W_1,\dots W_s$. Then
$$
\#\; \mbox{unfavorable cases}~ =~ \#\; \mbox{proper colorings of }\;\mathcal{G}_s~ \mbox{with  $m$ colors }~=~ \p_{\mathcal{G}_s}(m)
$$
where $\p_{\mathcal{G}_s}(m)$ is the chromatic polynomial of the graph $\mathcal{G}_s$.
Thus
$$
P_i(S)= 1- \frac{\p_{\mathcal{G}_s}(m)}{ m^w} \doteq q
$$
and therefore
$$
P(E_S)=q^N
$$
As before  two events
$E_S,E_{S'}\in {\cal S}$ are independent if $S\cap S' =\emptyset$. Therefore
The dependency graph for the family of  events ${\cal S}=\{E_S\}_{S\in P^*_{w}([n])}$
can be identified with the graph $G$ with vertex set  $ P^*_{w}([n])$ such that
two vertices  $S=\{W_1,\dots,W_s\}$ and $S'=\{W'_1,\dots,W'_s\}$  are connected by an edge of $G$ if
and only
$S\cap  S'\neq \emptyset$ (where recall that
 $S\cap S'\doteq
(\bigcup_{i=1}^s W_i)\cap  (\bigcup_{i=1}^s W'_i)$).
This implies that the neighbor $\Gamma^*_G(S)$ of a vertex  $S$ of $G$ is given by
$$
\Gamma_G^*(S)=\{{S'}: S'\in P^*_{w}([n])~{\rm and}~ S'\cap S\neq\emptyset\}
$$

\\By Theorem \ref{PAP}, the Moser-Tardos algorithm (as described in sec. \ref{halg}) finds  a separating hash family
${\rm SHF}(N; n,m,\{w_1,\dots,w_s\})$  if the
following  condition is satisfied:
there exists $\n>0$ such that
\be\label{FPs}
q^N\le \frac{\n}{  \displaystyle\sum_{\substack{Y\subseteq\Gamma^*_G(S)\\ Y~{\rm independent\,in\,} G}} \n^{|Y|}}
\ee
Note that,  as we did in the previous section,  we have set
$\m_S=\n$ for all  $S\in P^*_{w}([n])$.

\\The denominator of the r.h.s. of (\ref{FPs}) can be evaluated similarly as we did for the case of perfect hash families. Indeed,
given a disjoint family
$S=\{W_1,\dots,W_s\}$, the neighbor of $S$ in $G$ is formed by all vertices
$S'=\{W'_1,\dots,W'_s\}$ such that  $(\bigcup_{r=1}^s W_r)\bigcap (\,\bigcup_{r=1}^s W'_r)\neq \emptyset$.
The only thing that changes respect to the calculations done for case of the
perfect hash families is that now, fixed a set of columns $W$ with cardinality $w=w_1+\dots w_s$,  the number of different disjoint families
$S= \{W_1,\dots,W_s\}$
such that $W=\bigcup_{r=1}^s W_r$ and $|W_r|=w_r$ for $r=1,\dots, s$ is given by the quantity $m_w$ defined in (\ref{mw}).
Therefore we have
$$
\sum_{\substack{Y\subseteq\Gamma^*_G(S)\\ Y~{\rm independent\,in\,}G}} \n^{|Y|}=
1+ \sum_{k=1}^{\lfloor n/w\rfloor \wedge w} \G_k(w,n) (m_w\nu)^k
$$
where $\G_k(w,n)$ is exactly the same number defined in (\ref{gk}).
Hence posing $m_w\nu=\mu$ we have that  a   separating hash family ${\rm SHF}(N; n,m,\{w,\dots,w_s\})$  exists and can be found in polynomial time
by the Moser-Tardos algorithm if
\be\label{munu}
m_w q^N\le \frac{\m}{  \displaystyle1+ \sum_{k=1}^{\lfloor n/w\rfloor \wedge w} \G_k(w,n) ~\mu^k}
\ee
Note that the r.h.s. of inequality (\ref{munu}) and the r.h.s. of inequality (\ref{FP}) are the same.  Hence we get that the condition (\ref{FPs}) becomes
$$
\frac{(n-w)^{w-1}}{(w-1)!} m_w q^N\le \frac{1}{\varphi_{w,n}'(\t)}
$$
where $\varphi_{w,n}'(\t)$ is the same number as in (\ref{FP2}). In other word
the condition (\ref{FPs}) rewrites as
\be\label{SFP2}
N\ge \frac{ S_n(w)}{\ln\left(\frac{1}{  q}\right)}
\ee
with
\be\label{Snw}
S_n(w)= \ln[\varphi_{w,n}'(\t)]+(w-1)\ln{(n-w)}-\ln{(w-1)!}+ \ln (m_w)
\ee
According to Theorem \ref{t1} the MT-algorithm (as described in
section \ref{halg}) finds a ${\rm SHF}(N; n,m,\{w_1,\dots, w_s\})$ satisfying  (\ref{SFP2}) and hence (\ref{FPs}) in an
expected number of steps
$$\sum_{S\in P^*_{w}([n])} \n= |P^*_{w}([n])|\nu =
{n\choose w}m_w\nu={n\choose w}\m\le {n\choose w}$$
where the last inequality follows from the fact that
the optimum $\mu$
which maximize the r.h.s. of (\ref{munu}) is surely less than one. The number ${n\choose w}$ of these steps, similarly to what done
for PHF, has to be multiplied by $ N{w\choose 2}m_w{n\choose w}$ which is the number of operations needed to check item a) of step $i\ge 1$
of the MT-algorithm for SHF described in Sec. \ref{halg}.

\\This concludes the proof of Theorem \ref{shash}.

\section{Comparison with previous bounds}\label{compar}
\zeq
Let us
first
observe that the polynomial $\varphi_{w,n}(x)$ introduced in (\ref{fix}) is such that
$$
\lim_{n\to \infty}\varphi_{w,n}(x)\to (1+x)^w
$$
and therefore the number $\varphi'_{w,n}(\t)$ appearing in bounds (\ref{henne}) and (\ref{SHF}) for perfect
hash families and separating hash families resp.
 is such that
\be\label{asymf}
\lim_{n\to \infty} \varphi'_{w,n}(\t)= w\left(1+\frac{1}{ w-1}\right)^{w-1}
\ee
\subsection{Perfect Hash Families}
We first recall the previous lower bounds obtained in the literature via the Probabilistic method.
\\First, via the usual Lov\'asz local Lemma (see, e.g., \cite{DSW}) one obtains
\be\label{LLL}
N\ge \frac{L_n(w)}{ D_m(w)}
\ee
where
\be
L_n(w)=\ln\left[e\left(\binom{n}{ w}-\binom{n-w}{ w}\right)\right]\label{bn}
\ee
On the other hand, via the expurgation method (see \cite{STW} and \cite{DSW}) one gets
\be\label{expu}
N\ge \frac{E_n(w)}{ D_m(w)}
\ee
where
\be
E_n(w)=\ln \binom{2n}{ w}-\ln n
\ee
\\Let's first compare our bound (\ref{FP32}) with (\ref{expu}) obtained via expurgation method.
Note that the numerator  $E_n(w)$  appearing in the r.h.s. of
(\ref{expu}) can be written as
\begin{align}
E_n(w)&=\ln \binom{2n}{ w}-\ln n\nonumber\\
&= w\ln 2+ (w-1)\ln n+\sum_{j=1}^{w-1}\ln(1-\frac{j}{2n})-\ln(w!)\nonumber\\
&=w\ln 2+ (w-1)\ln (n-w)-\ln(w!)+\sum_{j=1}^{w-1}\ln\left(\frac{1-\frac{j}{ 2n}} {1-\frac{w}{ n}}\right)
\end{align}

\\So  that asymptotically
\be\label{asye}
E_n(w)~\sim~{\cal E}_w+ (w-1)\ln n-\ln[(w-1)!]
\ee
with
$$
{\cal E}_w= w\ln 2-\ln w
$$
On the other hand,
in view of (\ref{asymf}), the numerator of the r.h.s. of (\ref{FP32}) is asymptotic to
\be\label{asya}
A_n(w)~\sim ~ {\cal A}_w+(w-1)\ln{(n)}-\ln[(w-1)!]
\ee
with
$$
{\cal A}_w= \ln w+ \ln\left[\left(1+\frac{1}{ w-1}\right)^{w-1}\right]
$$
Observe that while ${\cal E}_w$ grows linearly in $w$, the factor ${\cal A}_w$ in (\ref{asya}) grows only logarithmically $w$.
Thus, to compare (\ref{FP32}) with (\ref{expu}) as $n\to \infty$  we have that
$$
A_n(w)<B_n(w)~~~~~~~~\Longleftrightarrow~~~~~~~~~ 2\ln w+ \ln\left[\left(1+\frac{1}{ w-1}\right)^{w-1}\right]~ <~w\ln 2
$$
which happens as soon as
$$
w~>~6.91043
$$
This means that, asymptotically in $n$,  for all $w\ge 7$ our bound beats the expurgation bound.
Moreover, numerical evidence suggests that the function
\be\label{deltan}
\D_n(w)=E_n(w)-A_n(w)= \ln\left(\frac{2^w}{w}\right) +\sum_{j=1}^{w-1}\ln\left(\frac{1-\frac{j}{ 2n}}{ 1-\frac{w}{ n}}\right)-\ln[\varphi_{w,n}'(\t)]
\ee
is decreasing in $n$ for  fixed $w$. If this were the case (we do not have a proof of that), our bound would be always better as long as  $w\ge 7$
(see Table 1).
For values of $w\leq 6$, one can perform numerical calculations with the bounds (\ref{FP32}) and (\ref{expu}) and see that
our bound  beats the bound obtained via expurgation only for low values of $n$
and the lower is $w$ the lower is the $n$ for which we win. In particular, we get that for $w=6$, $w=5$, $w=4$, and $w=3$
our bound beats expurgation bound
for all $n\le 97$, $n\le 34$,  $n\le 15$, and $n\le 7$ resp. (see also Table 2).

\begin{table}[ht]
\caption{Bounds on the cardinality of perfect hash families with $w\ge 7$} 
\centering 
\begin{tabular}{c c c c c c } 
\hline\hline 
$n$ & $m$ & $w$ & Theorem \ref{hash} & Expurgation \\ [0.5ex] 
\hline 
15 & 7 & 7 & 1437 & 1926 \\ 
50 & 7 & 7 & 3034 & 3191 \\
200 & 7 & 7 & 4529 & 4572 \\
1000& 7 & 7  & 6139 & 6152 \\
50 & 8 & 8 & 8463 & 9159 \\
200 & 8 & 8 & 12965 & 13282 \\
1000 & 8 & 8 & 17774 & 17988\\
1000 & 8 & 12 & 900 & 911\\
1000 & 8 & 50 & 53 & 54\\
1000 & 15 & 50 & 730 & 781\\
1000 & 18 & 50 & 2812 & 3037\\
[1ex] 
\hline 
\end{tabular}
\label{table1} 
\end{table}

\begin{table}[ht]
\caption{Bounds on the cardinality of perfect hash families with $w<7$} 
\centering 
\begin{tabular}{c c c c c c } 
\hline\hline 
$n$ & $m$ & $w$ & Theorem \ref{hash} & Expurgation \\ [0.5ex] 
\hline 
10 & 4 & 4 & 57 & 62 \\ 
15 & 4 & 4 & 76 & 77 \\
50 & 4 & 4 & 121 & 114 \\
10& 5 & 5 & 144 & 187 \\
15 & 5 & 5 & 211 & 234 \\
50 & 5 & 5 & 369 & 364 \\
15 & 6 & 6 & 558 & 681\\
50 & 6 & 6 & 1072 & 1092\\
90 & 6 & 6 & 1284 & 1287\\
200 & 6 & 6 & 1557 & 1546\\
[1ex] 
\hline 
\end{tabular}
\label{table2} 
\end{table}

\vskip.2cm
\\By theoretical reasons
(see \cite{BFPS}) the cluster expansion Lov\'asz Lemma is always better than the usual Lov\'asz Local Lemma,
so that bound (\ref{FP32}) always beats the bound (\ref{LLL}) for any pair $(w,n)$. In any case, it is  interesting to compare
our bound with the Lov\'asz Lemma bound asymptotically as $n\to \infty$.
We have that they are not equivalent, i.e. our bound beats LLL even asymptotically since
$$
\lim_{n\to \infty}\frac{\ln[\varphi_{n,w}'(\t)]+(w-1)\ln{(n-w)}-\ln{(w-1)!}}{\ln\left[e\left(\binom{n}{ w}-\binom{n-w}{ w}\right)\right]}=\frac{w-1}{ w}
$$

\subsection{Separating Hash Families}
We can compare the bound (\ref{SHF}) obtained in Theorem \ref{shash}
with that obtained via expurgation method and see that the improvement is completely analogous
to that obtained for the perfect hash families.
First observe that, recalling (\ref{mw}), the quantity $S_n(w)$ defined in (\ref{Snw}) can be rewritten as follows
$$
S_n(w)
=  \ln[w\varphi_{w,n}'(\t)] +(w-1)\ln{(n-w)}- \sum_{i=1}^s\ln(w_i!) - \sum_{p=1}^w\ln(m_p!)
$$
Let us then consider for simplicity the case $k=2$ and $w_1\neq w_2$. For this case
 $S_n(w)$  becomes
$$
S_n(w)=\ln[w\varphi_{w,n}'(\t)] +(w-1)\ln{(n-w)}- \ln(w_1!w_2!)
$$
Via the expurgation method (see  \cite{DSW}) one gets
\be\label{expu2}
N\ge\frac{F_n(w)}{ \ln\left(\frac{1}{  q}\right)}
\ee
where
\be
F_n(w)=\ln \binom{2n}{w_1}-\ln \binom{2n-w_1}{w_2}-\ln n
\ee
As before we can write
$$
F_n(w)=w\ln 2+ (w-1)\ln (n-w)+\sum_{j=1}^{w-1}\ln\left(\frac{1-\frac{j}{ 2n}}{ 1-\frac{w}{ n}}\right)-\ln(w_1!w_2!)
$$
So we have
$$
F_n(w) - S_n(w)=\D_n(w)
$$
where $\D_n(w)$ is  the same quantity defined in (\ref{deltan}). This means that all that we discussed for perfect hash families
holds also for separating hash families. In particular, our bound beats the bound obtained via expurgation method  reported in \cite{DSW} asymptotically in $n$ as soon as
$w>6$.

\\We finally compare our bound with still another bound  given by Stinson and Zaverucha \cite{SZ} in 2008. These authors
claim that  a 
${\rm SHF}(N; n,m,\{w_1,w_2\})$  exists if
\be\label{SZ}
n\le \left( 1- \frac{1}{  C_w}\right)\left[\frac{1}{ q}\right]^\frac{N}{ w-1}
\ee
where $w=w_1+w_2$ and
$$
C_w=\begin{cases} w_1!w_2! & \text{if } w_1\neq w_2\\
 2 w_1!w_2!  & \text{if } w_1= w_2
 \end{cases}
$$
Our bound (\ref{SHF}) on the other hand implies that  a 
$ {\rm SHF}(N; n,m,\{w_1,w_2\})$  exists if
$$
(n-w)^{w-1} \le \frac{{C_w}}{ w\varphi'(\t)} \left[\frac{1}{ q}\right]^{N}
$$
Once again, for sake of simplicity we perform this comparison asymptotically as $N$ (and hence $n$) large. In this case we have seen that
$ \varphi'(\t)\sim w(1+\frac{1}{ w-1})^{w-1}$. Therefore asymptotically  we have the bound
\be\label{vai}
n\le \left(\frac{C_w}{ w^2}\right)^{\frac{1}{ w-1}}\frac{w-1}{ w} \left[\frac{1}{ q}\right]^{\frac{N}{ w-1}}
\ee
The comparison of this with bound (\ref{SZ}) in now straightforward. We see that our bound (\ref{vai})  beats  (\ref{SZ}) as soon as
$w$ is larger than $6$.


%
%
%
%
%



\section*{Acknowledgments}

\\ A.P.  has been partially supported by the Brazilian  agencies
Conselho Nacional de Desenvolvimento Cient\'{\i}fico e Tecnol\'ogico
(CNPq  - Bolsa de Produtividade em pesquisa, grant n. 306208/2014-8)
and  Funda{\c{c}}\~ao de Amparo \`a  Pesquisa do Estado de Minas Gerais (FAPEMIG - Programa de Pesquisador Mineiro, grant n. 00230/14).
R. S. has been partially supported by
Conselho Nacional de Desenvolvimento Cient\'{\i}fico e Tecnol\'ogico
(CNPq  - Bolsa de Produtividade em pesquisa, grant n. 310552/2014-1). The authors thank
Cristiano Santos Benjamin for his valuable  help
with the numerical computations present in the paper.

%

\end{document}